\documentclass[12pt]{article}%
\usepackage{amsfonts}
\usepackage{sw20bams}
\usepackage{amsmath}
\usepackage{amssymb}
\usepackage{graphicx}%
\providecommand{\U}[1]{\protect\rule{.1in}{.1in}}
\begin{document}

\title{Rounds, Color, Parity, Squares}
\author{Steven Finch}
\date{January 23, 2022}
\maketitle

\begin{abstract}
This is a sequel to our paper \textquotedblleft Permute, Graph, Map,
Derange\textquotedblright,\ involving decomposable combinatorial labeled
structures in the exp-log class of type $a=1/2$, $1$, $3/2$, $2$. As before,
our approach is to establish how well existing theory matches experimental
data and to raise open questions.

\end{abstract}

\footnotetext{Copyright \copyright \ 2021--2022 by Steven R. Finch. All rights
reserved.}We continue where we left off in \cite{Fi-tcs7}. \ Among the most
striking features of a combinatorial object with $n$~nodes are

\begin{itemize}
\item the number of cycles or components,

\item the size of the longest cycle or largest component,

\item the size of the shortest cycle or smallest component.
\end{itemize}

\noindent The latter two topics will be our focus. \ Throughout this paper, a
random object is chosen uniformly from a set (to be explicated in the
following sections). \ Let $b_{n}$ denote the number of $n$-objects and
$c_{n}$ denote the number of $n$-objects that are connected, i.e., who possess
exactly one component.

Key to our study are recursive formulas for $L_{k,n}$ and $S_{k,n}$, the
number of $n$-objects whose largest and smallest components, respectively,
have exactly $k$ nodes, $1\leq k\leq n$. \ The initial conditions are%
\[%
\begin{array}
[c]{ccc}%
L_{0,n}=\delta_{0,n}, &  & L_{1,n}=\left(  1-\delta_{0,n}\right)  c_{1}^{n}%
\end{array}
\]
and%
\[%
\begin{array}
[c]{c}%
S_{0,n}=\delta_{0,n}c_{1},
\end{array}
\]
respectively. \ Letting
\[
m_{j,k,n}=\min\{k-1,n-k\,j\}
\]
and suppressing dependence on $j$, $k$, $n$, we have \cite{Fi-tcs7, GG-tcs7,
PR1-tcs7}%
\[
L_{k,n}=%
{\displaystyle\sum\limits_{j=1}^{\left\lfloor n/k\right\rfloor }}
\frac{n!c_{k}^{j}}{j!(k!)^{j}(n-k\,j)!}\cdot\left\{
\begin{array}
[c]{ccc}%
{\displaystyle\sum\limits_{i=1}^{m}}
L_{i,n-k\,j} &  & \text{if }m\geq1,\\%
{\displaystyle\sum\limits_{i=0}^{1}}
L_{i,n-k\,j} &  & \text{if }m=0.
\end{array}
\right.
\]
Letting%
\[
\theta_{k,n}=\left\{
\begin{array}
[c]{lll}%
1 &  & \text{if }k\text{ is a divisor of }n\text{,}\\
0 &  & \text{otherwise}%
\end{array}
\right.
\]
we have \cite{PR1-tcs7}%
\[
S_{k,n}=%
{\displaystyle\sum\limits_{j=1}^{\left\lfloor n/k\right\rfloor }}
\frac{n!c_{k}^{j}}{j!(k!)^{j}(n-k\,j)!}\cdot%
{\displaystyle\sum\limits_{i=k+1}^{n-k\,j}}
S_{i,n-k\,j}+\theta_{k,n}\frac{n!c_{k}^{n/k}}{(n/k)!(k!)^{n/k}}.
\]
Clearly
\[%
{\displaystyle\sum\limits_{k=1}^{n}}
L_{k,n}=%
{\displaystyle\sum\limits_{k=1}^{n}}
S_{k,n}=b_{n}%
\]
and $L_{n,n}=S_{n,n}=c_{n}$. \ A computer algebra software package (e.g.,
Mathematica) makes exact integer calculations for ample $n$ of $L_{k,n}$ and
$S_{k,n}$ feasible.

Permutations and derangements belong to the exp-log class of type $a=1$,
whereas mappings belong to the exp-log class of type $a=1/2$. \ Explaining the
significance of the parameter $a>0$ would take us too far afield
\cite{PR2-tcs7}. \ Let%
\[%
\begin{array}
[c]{ccc}%
E(x)=%
{\displaystyle\int\limits_{x}^{\infty}}
\dfrac{e^{-t}}{t}dt=-\operatorname{Ei}(-x), &  & x>0
\end{array}
\]
be the exponential integral. \ Define \cite{SL-tcs7, Shi-tcs7, FO-tcs7,
Gou-tcs7, ABT-tcs7, Pin-tcs7}%
\[
_{L}G_{a}(r,h)=\frac{\Gamma(a+1)a^{r-1}}{\Gamma(a+h)(r-1)!}%
{\displaystyle\int\limits_{0}^{\infty}}
x^{h-1}E(x)^{r-1}\exp\left[  -a\,E(x)-x\right]  dx,
\]%
\[
_{S}G_{a}(r,h)=\left\{
\begin{array}
[c]{lll}%
e^{-h\,\gamma}a^{r-1}/r! &  & \text{if }h=a,\\
\dfrac{\Gamma(a+1)}{(h-1)!(r-1)!}%
{\displaystyle\int\limits_{0}^{\infty}}
x^{h-1}\exp\left[  a\,E(x)-x\right]  dx &  & \text{if }h>a
\end{array}
\right.
\]
which are related to the $h^{\text{th}}$ moment of the $r^{\text{th}}$
largest/smallest component size (in this paper, rank $r=1$ or $2$; height
$h=1$ or $2$). \ Our notation $_{S}G_{a}$ is deceiving. \ While permutation
and\ derangement moments coincide for $L$ (both being $_{L}G_{a}$ with $a=1$),
they are \textbf{not} equal for $S$ (they differ by a factor $e$). \ 

We need a certain equation. \ For $a>0$, the asymptotic probability that the
largest component has size $>n\,x$ is \cite{Ew-tcs7}%
\[
a%
{\displaystyle\int\limits_{x}^{1}}
\frac{(1-y)^{a-1}}{y}dy=\frac{1}{2},
\]
a condition which is met when%
\[
x=\left\{
\begin{array}
[c]{lll}%
4e/(1+e)^{2} &  & \text{if }a=1/2,\\
1/\sqrt{e} &  & \text{if }a=1,\\
\operatorname{sech}(\xi)^{2} &  & \text{if }a=3/2,\\
-W\left(  -e^{-5/4}\right)  &  & \text{if }a=2;
\end{array}
\right.
\]
$\xi$ is the unique positive solution of $\tanh(\xi)=-1/6+\xi$ and $W$ is
Lambert's function.

For fixed $n$, the sequences $\{L_{k,n}:1\leq k\leq n\}$ and $\{S_{k,n}:1\leq
k\leq n\}$ constitute probability mass functions (upon normalization by
$b_{n}$). \ These have corresponding means $_{L}\mu_{n}$, $_{S}\mu_{n}$ and
variances $_{L}\sigma_{n}^{2}$, $_{S}\sigma_{n}^{2}$ given in the tables. \ We
also provide the median $_{L}\nu_{n}$; note that $_{S}\nu_{n}=1$ for $n>5$ is
trivial. \ For convenience (in table headings only), the following notation is
used:
\[%
\begin{array}
[c]{ccccc}%
_{L}\widetilde{\mu}_{n}=\dfrac{_{L}\mu_{n}}{n}, &  & _{L}\widetilde{\sigma
}_{n}^{2}=\dfrac{_{L}\sigma_{n}^{2}}{n^{2}}, &  & _{L}\widetilde{\nu}%
_{n}=\dfrac{_{L}\nu_{n}}{n},
\end{array}
\]%
\[%
\begin{array}
[c]{ccc}%
_{S}\widetilde{\mu}_{n}=\left\{
\begin{array}
[c]{lll}%
\dfrac{_{S}\mu_{n}}{n^{1/2}} & \bigskip & \text{if }a=1/2,\\
\dfrac{_{S}\mu_{n}}{\ln(n)} & \bigskip & \text{if }a=1,\\
\dfrac{_{S}\mu_{n}}{1} & \bigskip & \text{if }a=3/2\text{ or }2;
\end{array}
\right.  &  & _{S}\widetilde{\sigma}_{n}^{2}=\left\{
\begin{array}
[c]{lll}%
\dfrac{_{S}\sigma_{n}^{2}}{n^{3/2}} & \bigskip & \text{if }a=1/2,\\
\dfrac{_{S}\sigma_{n}^{2}}{n} & \bigskip & \text{if }a=1,\\
\dfrac{_{S}\sigma_{n}^{2}}{n^{1/2}} & \bigskip & \text{if }a=3/2,\\
\dfrac{_{S}\sigma_{n}^{2}}{\ln(n)} &  & \text{if }a=2.
\end{array}
\right.
\end{array}
\]
Finally, we wonder about the existence of \textquotedblleft
natural\textquotedblright\ combinatorial objects in exp-log class of type
$a=1/3$ or $5/4$, say. \ These would be exceedingly interesting to examine.

\section{Children's Rounds}

Consider the myriad arrangements of $n$ labeled children into rounds
\cite{St1-tcs7, St2-tcs7, O1-tcs7}. \ A \textit{round} means the same as a
directed ring or circle, with exactly one child inside the ring and the others
encircling. \ We permit the outer ring to have as few as one child. \ For
$n=4$, there are $8$ ways to make one round (a $4$-round of children, $1$
inside and $3$ outside):%
\[%
\begin{array}
[c]{ccccccc}%
\begin{array}
[c]{ccc}
& 4 & \\
& 1 & \\
2 &  & 3
\end{array}
&  &
\begin{array}
[c]{ccc}
& 3 & \\
& 1 & \\
2 &  & 4
\end{array}
&  &
\begin{array}
[c]{ccc}
& 4 & \\
& 2 & \\
1 &  & 3
\end{array}
&  &
\begin{array}
[c]{ccc}
& 3 & \\
& 2 & \\
1 &  & 4
\end{array}
\end{array}
\]%
\[%
\begin{array}
[c]{ccccccc}%
\begin{array}
[c]{ccc}
& 4 & \\
& 3 & \\
1 &  & 2
\end{array}
&  &
\begin{array}
[c]{ccc}
& 2 & \\
& 3 & \\
1 &  & 4
\end{array}
&  &
\begin{array}
[c]{ccc}
& 3 & \\
& 4 & \\
1 &  & 2
\end{array}
&  &
\begin{array}
[c]{ccc}
& 2 & \\
& 4 & \\
1 &  & 3
\end{array}
\end{array}
\]
and $12$ ways to make two rounds (both $2$-rounds):%
\[%
\begin{array}
[c]{ccccccccccc}%
\begin{array}
[c]{ccc}%
1 & | & 3\\
2 & | & 4
\end{array}
&  &
\begin{array}
[c]{ccc}%
1 & | & 4\\
2 & | & 3
\end{array}
&  &
\begin{array}
[c]{ccc}%
1 & | & 2\\
3 & | & 4
\end{array}
&  &
\begin{array}
[c]{ccc}%
1 & | & 4\\
3 & | & 2
\end{array}
&  &
\begin{array}
[c]{ccc}%
1 & | & 2\\
4 & | & 3
\end{array}
&  &
\begin{array}
[c]{ccc}%
1 & | & 3\\
4 & | & 2
\end{array}
\end{array}
\]%
\[%
\begin{array}
[c]{ccccccccccc}%
\begin{array}
[c]{ccc}%
2 & | & 3\\
1 & | & 4
\end{array}
&  &
\begin{array}
[c]{ccc}%
2 & | & 4\\
1 & | & 3
\end{array}
&  &
\begin{array}
[c]{ccc}%
3 & | & 2\\
1 & | & 4
\end{array}
&  &
\begin{array}
[c]{ccc}%
3 & | & 4\\
1 & | & 2
\end{array}
&  &
\begin{array}
[c]{ccc}%
4 & | & 2\\
1 & | & 3
\end{array}
&  &
\begin{array}
[c]{ccc}%
4 & | & 3\\
1 & | & 2
\end{array}
\end{array}
\]
implying that $L_{2,4}=12=S_{2,4}$ and $L_{4,4}=8=S_{4,4}$. \ For $n=5$, there
are $30$ ways to make a $5$-round and $60$ ways to make a $3$-round \& a
$2$-round, implying that $L_{3,5}=60=S_{2,5}$ and $L_{5,5}=30=S_{5,5}$. \ For
$n=6$, there are $144$ ways to make a $6$-round, $240$ ways to make a
$4$-round \& a $2$-round, $90$ ways to make two $3$-rounds, and $120$ ways to
make three $2$-rounds, implying that $L_{2,6}=120$, $S_{2,6}=360$,
$L_{3,6}=90=S_{3,6}$, $L_{4,6}=240$ and $L_{6,6}=144=S_{6,6}$. \ We obtain%
\[%
\begin{array}
[c]{ccc}%
c_{n}=(n-1)!+(n-2)!, &  & n\geq2
\end{array}
\]
and, upon normalization by $b_{n}$,

\begin{center}%
\begin{tabular}
[c]{|c|c|c|c|c|c|}\hline
$n$ & $_{L}\widetilde{\mu}_{n}$ & $_{L}\widetilde{\sigma}_{n}^{2}$ &
$_{L}\widetilde{\nu}_{n}$ & $_{S}\widetilde{\mu}_{n}$ & $_{S}\widetilde
{\sigma}_{n}^{2}$\\\hline
1000 & 0.621184 & 0.036672 & 0.6020 & 0.862134 & 1.317448\\\hline
2000 & 0.622539 & 0.036764 & 0.6040 & 0.834706 & 1.312715\\\hline
3000 & 0.623052 & 0.036802 & 0.6050 & 0.820866 & 1.311031\\\hline
4000 & 0.623326 & 0.036823 & 0.6053 & 0.811867 & 1.310156\\\hline
\end{tabular}

Table 1A:\ Statistics for Children's Rounds ($a=1$)
\end{center}

\noindent also%
\[
\lim_{n\rightarrow\infty}\dfrac{_{L}\mu_{n}}{n}=\,_{L}G_{1}%
(1,1)=0.62432998854355087099...,
\]%
\[
\lim_{n\rightarrow\infty}\dfrac{_{L}\sigma_{n}^{2}}{n^{2}}=\,_{L}%
G_{1}(1,2)-\,_{L}G_{1}(1,1)^{2}=0.03690783006485220217...,
\]%
\[
\lim_{n\rightarrow\infty}\dfrac{_{L}\nu_{n}}{n}=\frac{1}{\sqrt{e}%
}=0.60653065971263342360...,
\]%
\[
\lim_{n\rightarrow\infty}\dfrac{_{S}\mu_{n}}{\ln(n)}=e^{-\gamma}%
=0.56145948356688516982...,
\]%
\[
\lim_{n\rightarrow\infty}\dfrac{_{S}\sigma_{n}^{2}}{n}=\,_{S}G_{1}%
(1,2)=1.30720779891056809974....
\]

\noindent It is not surprising that $_{S}\sigma_{n}^{2}$\ enjoys linear
growth:\ $S_{1,n}\sim(1-1/e)n!$ and $S_{n,n}=(n-1)!$ jointly place
considerable weight on the distributional extremes. \ The unusual logarithmic
growth of $_{S}\mu_{n}$\ is due to $S_{1,n}$ nevertheless overwhelming all
other $S_{k,n}$.

\subsection{Variant}

If we disallow the outer ring from having just one child, then clearly
$L_{2,n}$, $S_{2,n}$ and $L_{n-2,n}$ are all zero. We obtain%
\[
c_{n}=\left\{
\begin{array}
[c]{lll}%
0 &  & \text{if }n=2,\\
(n-1)!+(n-2)! &  & \text{if }n\geq3
\end{array}
\right.
\]

\begin{center}%
\begin{tabular}
[c]{|c|c|c|c|c|c|}\hline
$n$ & $_{L}\widetilde{\mu}_{n}$ & $_{L}\widetilde{\sigma}_{n}^{2}$ &
$_{L}\widetilde{\nu}_{n}$ & $_{S}\widetilde{\mu}_{n}$ & $_{S}\widetilde
{\sigma}_{n}^{2}$\\\hline
1000 & 0.622431 & 0.036818 & 0.6030 & 1.846026 & 3.574315\\\hline
2000 & 0.623163 & 0.036837 & 0.6050 & 1.816840 & 3.564892\\\hline
3000 & 0.623468 & 0.036851 & 0.6053 & 1.802115 & 3.561459\\\hline
4000 & 0.623638 & 0.036860 & 0.6055 & 1.792542 & 3.559653\\\hline
\end{tabular}

Table 1B:\ Statistics for Variant of Children's Rounds ($a=1$)
\end{center}

\noindent and, while the $L$ limits are the same as before, the $S$ limits
differ by a factor of $e$:%
\[
\lim_{n\rightarrow\infty}\dfrac{_{S}\mu_{n}}{\ln(n)}=e^{1-\gamma
}=1.52620511159586388047...,
\]%
\[
\lim_{n\rightarrow\infty}\dfrac{_{S}\sigma_{n}^{2}}{n}=e\cdot\,_{S}%
G_{1}(1,2)=3.55335920579854297440...
\]
as argued in Section 4 of \cite{Fi-tcs7}.

\section{Cycle-Colored Permutations}

Assuming two colors are available \cite{O2-tcs7}, it is clear that%
\[%
\begin{array}
[c]{ccc}%
c_{n}=2\cdot(n-1)!, &  & n\geq1.
\end{array}
\]
Let us explain why $L_{2,3}=12$ and $S_{1,3}=20$. \ The $3$ permutations
\[%
\begin{array}
[c]{ccccc}%
(1\;2)(3), &  & (1\;3)(2) &  & (2\;3)(1)
\end{array}
\]
each have $4$ possible colorings, giving $12$ to both $L_{2,3}$ and $S_{1,3}$;
\ the permutation $(1)(2)(3)$ contributes another $8$ to $S_{1,3}$. In
contrast, $L_{2,4}=60$, $L_{3,4}=32$ and $S_{2,4}=12$ because the $3$
permutations
\[%
\begin{array}
[c]{ccccc}%
(1\;2)(3\;4), &  & (1\;3)(2\;4) &  & (1\;4)(2\;3)
\end{array}
\]
again give $12$ to $L_{2,4}$ and $S_{2,4}$, whereas the $8$ permutations%
\[%
\begin{array}
[c]{ccccccc}%
(1\;2\;3)(4), &  & (1\;3\;2)(4), &  & (1\;2\;4)(3), &  & (1\;4\;2)(3),
\end{array}
\]%
\[%
\begin{array}
[c]{ccccccc}%
(1\;3\;4)(2), &  & (1\;4\;3)(2), &  & (2\;3\;4)(1), &  & (2\;4\;3)(1)
\end{array}
\]
give $32$ to $L_{3,4}$; also, the $6$ permutations%
\[%
\begin{array}
[c]{ccccc}%
(1)(2)(3\;4), &  & (1)(3)(2\;4) &  & (1)(4)(2\;3)
\end{array}
\]%
\[%
\begin{array}
[c]{ccccc}%
(3)(4)(1\;2), &  & (2)(4)(1\;3) &  & (2)(3)(1\;4)
\end{array}
\]
contribute another $48$ to $L_{2,4}$.

Upon normalization by $(n+1)!$, we obtain

\begin{center}%
\begin{tabular}
[c]{|c|c|c|c|c|c|}\hline
$n$ & $_{L}\widetilde{\mu}_{n}$ & $_{L}\widetilde{\sigma}_{n}^{2}$ &
$_{L}\widetilde{\nu}_{n}$ & $_{S}\widetilde{\mu}_{n}$ & $_{S}\widetilde
{\sigma}_{n}^{2}$\\\hline
1000 & 0.476115 & 0.027160 & 0.4480 & 1.292899 & 1.011228\\\hline
2000 & 0.475877 & 0.027132 & 0.4480 & 1.291149 & 0.976960\\\hline
3000 & 0.475798 & 0.027123 & 0.4480 & 1.290504 & 0.959564\\\hline
4000 & 0.475758 & 0.027119 & 0.4480 & 1.290163 & 0.948225\\\hline
\end{tabular}

Table 2A:\ Statistics for Cycle-Colored Permutations ($a=2$)
\end{center}

\noindent and%

\[
\lim_{n\rightarrow\infty}\dfrac{_{L}\mu_{n}}{n}=\,_{L}G_{2}%
(1,1)=0.47563939666525000670...,
\]%
\[
\lim_{n\rightarrow\infty}\dfrac{_{L}\sigma_{n}^{2}}{n^{2}}=\,_{L}%
G_{2}(1,2)-\,_{L}G_{2}(1,1)^{2}=0.02710536578919830440...,
\]%
\[
\lim_{n\rightarrow\infty}\dfrac{_{L}\nu_{n}}{n}=-W\left(  -e^{-5/4}\right)
=0.44878202648462460223...,
\]%
\[
\lim_{n\rightarrow\infty}\dfrac{_{S}\mu_{n}}{1}=\kappa=1.29...,
\]%
\[
\lim_{n\rightarrow\infty}\dfrac{_{S}\sigma_{n}^{2}}{\ln(n)}=e^{-2\gamma
}=0.31523675168719339806....
\]

\noindent No explicit integral is known for the constant $\kappa$ -- it
appears again shortly -- the $e^{-2\gamma}$ limit is proved in Theorem 5 of
\cite{PR2-tcs7}.

\subsection{Variant}

If we prohibit $1$-cycles, then clearly $L_{1,n}$, $S_{1,n}$ and $L_{n-1,n}$
are all zero; however $L_{2,4}=12=S_{2,4}.$ \ We obtain%
\[
c_{n}=\left\{
\begin{array}
[c]{lll}%
0 &  & \text{if }n=1,\\
2(n-1)! &  & \text{if }n\geq2
\end{array}
\right.
\]

\begin{center}%
\begin{tabular}
[c]{|c|c|c|c|c|c|}\hline
$n$ & $_{L}\widetilde{\mu}_{n}$ & $_{L}\widetilde{\sigma}_{n}^{2}$ &
$_{L}\widetilde{\nu}_{n}$ & $_{S}\widetilde{\mu}_{n}$ & $_{S}\widetilde
{\sigma}_{n}^{2}$\\\hline
1000 & 0.477065 & 0.027268 & 0.4480 & 3.159931 & 6.534345\\\hline
2000 & 0.476353 & 0.027187 & 0.4485 & 3.149165 & 6.372009\\\hline
3000 & 0.476115 & 0.027160 & 0.4483 & 3.145123 & 6.288169\\\hline
4000 & 0.475996 & 0.027146 & 0.4482 & 3.142963 & 6.233112\\\hline
\end{tabular}

Table 2B:\ Statistics for Cycle-Colored Derangements ($a=2$)
\end{center}

\noindent and, while the $L$ limits are the same as before, the $S$ limits
differ somewhat:%
\[
\lim_{n\rightarrow\infty}\dfrac{_{S}\mu_{n}}{1}=\frac{\kappa-1+e^{-2}}{e^{-2}%
}=3.13...,
\]%
\[
\lim_{n\rightarrow\infty}\dfrac{_{S}\sigma_{n}^{2}}{\ln(n)}=e^{2-2\gamma
}=2.32930204266134332103....
\]
The expresson involving $\kappa$ for the average shortest cycle length follows
from%
\[
\frac{_{S}\mu_{n}\cdot b_{n}+1\cdot\left(  (n+1)!-b_{n}\right)  }{(n+1)!}%
\sim\kappa
\]
and the fact that $b_{n}/(n+1)!\rightarrow1/e^{2}$ as $n\rightarrow\infty$;
hence
\[
_{S}\mu_{n}\cdot e^{-2}+\left(  1-e^{-2}\right)  \sim\kappa.
\]
Similarly,%
\[
\frac{_{S}\sigma_{n}^{2}\cdot b_{n}+1\cdot\left(  (n+1)!-b_{n}\right)
}{(n+1)!}\sim e^{-2\gamma}\ln(n)
\]
yields%
\[
_{S}\sigma_{n}^{2}\cdot e^{-2}\sim e^{-2\gamma}\ln(n),
\]
although the rightmost column of Table 2B suggests that this approximation is poor.

\section{Component-Colored Mappings}

Assuming three colors are available \cite{O5-tcs7}, it is clear that%
\[%
\begin{array}
[c]{ccc}%
c_{n}=3\cdot n!%
{\displaystyle\sum\limits_{j=1}^{n}}
\dfrac{n^{n-j-1}}{(n-j)!}, &  & n\geq1.
\end{array}
\]
An argument in Section 3 of \cite{Fi-tcs7} gives $L_{1,3}=27$, $L_{2,3}=81$,
$S_{1,3}=108$ and $L_{3,3}=51=S_{3,3}$. \ Similarly, $L_{1,4}=81$,
$L_{2,4}=729$, $L_{3,4}=612$, $S_{1,4}=1179$, $S_{2,4}=243$ and $L_{4,4}%
=426=S_{4,4}$. \ 

Upon normalization by $b_{n}$, we obtain

\begin{center}%
\begin{tabular}
[c]{|c|c|c|c|c|c|}\hline
$n$ & $_{L}\widetilde{\mu}_{n}$ & $_{L}\widetilde{\sigma}_{n}^{2}$ &
$_{L}\widetilde{\nu}_{n}$ & $_{S}\widetilde{\mu}_{n}$ & $_{S}\widetilde
{\sigma}_{n}^{2}$\\\hline
1000 & 0.544944 & 0.032583 & 0.5170 & 2.590160 & 5.925381\\\hline
2000 & 0.542744 & 0.032407 & 0.5155 & 2.597466 & 6.079830\\\hline
3000 & 0.541778 & 0.032331 & 0.5150 & 2.600326 & 6.152861\\\hline
4000 & 0.541205 & 0.032285 & 0.5142 & 2.601910 & 6.198088\\\hline
\end{tabular}

Table 3:\ Statistics for Component-Colored Mappings ($a=3/2$)
\end{center}

\noindent and%
\[
\lim_{n\rightarrow\infty}\dfrac{_{L}\mu_{n}}{n}=\,_{L}G_{3/2}%
(1,1)=0.53753956799272857702...,
\]%
\[
\lim_{n\rightarrow\infty}\dfrac{_{L}\sigma_{n}^{2}}{n^{2}}=\,_{L}%
G_{3/2}(1,2)-\,_{L}G_{3/2}(1,1)^{2}=0.0319941833398955610...,
\]%
\[
\lim_{n\rightarrow\infty}\dfrac{_{L}\nu_{n}}{n}=\operatorname{sech}(\xi
)^{2}=0.51092489978120431153...,
\]%
\[
\lim_{n\rightarrow\infty}\dfrac{_{S}\mu_{n}}{1}=2.61...,
\]%
\[
\lim_{n\rightarrow\infty}\dfrac{_{S}\sigma_{n}^{2}}{n^{1/2}}=6.50....
\]
No explicit integrals are known for the latter two results.

\subsection{Variant}

Instead of removing the smallest possible components ($1$-components), we
wonder about removing the largest possible components. \ Of course, this can
be done only imprecisely, as the size of the giant component cannot be known
beforehand. \ If a \textquotedblleft post-processing removal\textquotedblright%
\ is acceptable (rather than a \textquotedblleft pre-processing
removal\textquotedblright), then for the remaining components\ we easily have
\cite{Gou-tcs7}
\[
\lim_{n\rightarrow\infty}\dfrac{_{L}\mu_{n}}{n}=\,_{L}G_{3/2}%
(2,1)=0.21627840128093004373...,
\]%
\[
\lim_{n\rightarrow\infty}\dfrac{_{L}\sigma_{n}^{2}}{n^{2}}=\,_{L}%
G_{3/2}(2,2)-\,_{L}G_{3/2}(2,1)^{2}=0.00867133690820287157....
\]
Both mean and variance are significantly reduced (from $0.537$ to $0.216$ and
$0.032$ to $0.009$, respectively).

The recursive formulas for $L_{k,n}$ and $S_{k,n}$, however, apply expressly
when the desired rank $r=1$. \ We cannot use our current exact integer-based
algorithm to experimentally confirm these statistics for $r=2$. \ A\ Monte
Carlo simulation would be feasible, but less thorough and less accurate. \ 

\section{Parity of Cycle Lengths}

Let EV\ and OD\ refer to permutations with all cycle lengths even and with all
cycle lengths odd, respectively \cite{Bn-tcs7, Lg-tcs7, O3-tcs7}. \ Clearly%
\[
c_{n}=\left\{
\begin{array}
[c]{lll}%
(n-1)! &  & \text{if }n\equiv0\operatorname{mod}2,\\
0 &  & \text{if otherwise}%
\end{array}
\right.
\]
holds for EV permutations and
\[
c_{n}=\left\{
\begin{array}
[c]{lll}%
(n-1)! &  & \text{if }n\equiv1\operatorname{mod}2,\\
0 &  & \text{if otherwise}%
\end{array}
\right.
\]
holds for OD\ permutations. \ On the one hand, in EV we have $L_{k,n}%
=0=S_{k,n}$ for all $k$ and\ odd $n$. On the other hand, in OD\ we have
$L_{1,n}=1$ always (due to the identity permutation) and $L_{3,4}=8$,
$S_{3,6}=40$ (both quickly checked). \ The fact that $S_{n,n}=0$ for even $n$
but $S_{n,n}=(n-1)!$ for odd $n$ gives the divergence between Tables 4B \&\ 4C
for OD.

Upon normalization by $b_{n}$, we obtain

\begin{center}%
\begin{tabular}
[c]{|c|c|c|c|c|c|}\hline
$n$ & $_{L}\widetilde{\mu}_{n}$ & $_{L}\widetilde{\sigma}_{n}^{2}$ &
$_{L}\widetilde{\nu}_{n}$ & $_{S}\widetilde{\mu}_{n}$ & $_{S}\widetilde
{\sigma}_{n}^{2}$\\\hline
1000 & 0.758202 & 0.037044 & 0.7850 & 2.028405 & 1.400424\\\hline
2000 & 0.758012 & 0.037026 & 0.7865 & 2.037816 & 1.400250\\\hline
3000 & 0.757949 & 0.037020 & 0.7863 & 2.042013 & 1.400192\\\hline
4000 & 0.757918 & 0.037016 & 0.7862 & 2.044523 & 1.400163\\\hline
\end{tabular}

Table 4A:\ Statistics for EV permutations ($a=1/2$)
\end{center}

\noindent

\begin{center}%
\begin{tabular}
[c]{|c|c|c|c|c|c|}\hline
$n$ & $_{L}\widetilde{\mu}_{n}$ & $_{L}\widetilde{\sigma}_{n}^{2}$ &
$_{L}\widetilde{\nu}_{n}$ & $_{S}\widetilde{\mu}_{n}$ & $_{S}\widetilde
{\sigma}_{n}^{2}$\\\hline
1000 & 0.757601 & 0.036937 & 0.7860 & 0.552117 & 0.125460\\\hline
2000 & 0.757712 & 0.036972 & 0.7860 & 0.553094 & 0.125440\\\hline
3000 & 0.757749 & 0.036984 & 0.7860 & 0.553535 & 0.125434\\\hline
4000 & 0.757768 & 0.036990 & 0.7860 & 0.553800 & 0.125431\\\hline
\end{tabular}

Table 4B:\ Statistics for OD permutations ($a=1/2$, $n\equiv
0\operatorname{mod}2$)
\end{center}

\noindent

\begin{center}%
\begin{tabular}
[c]{|c|c|c|c|c|c|}\hline
$n$ & $_{L}\widetilde{\mu}_{n}$ & $_{L}\widetilde{\sigma}_{n}^{2}$ &
$_{L}\widetilde{\nu}_{n}$ & $_{S}\widetilde{\mu}_{n}$ & $_{S}\widetilde
{\sigma}_{n}^{2}$\\\hline
999 & 0.758045 & 0.037077 & 0.7845 & 1.501395 & 1.274342\\\hline
1999 & 0.757934 & 0.037042 & 0.7864 & 1.502628 & 1.274497\\\hline
2999 & 0.757897 & 0.037031 & 0.7863 & 1.503154 & 1.274549\\\hline
3999 & 0.757878 & 0.037025 & 0.7862 & 1.503461 & 1.274575\\\hline
\end{tabular}

Table 4C:\ Statistics for OD permutations ($a=1/2$, $n\equiv
1\operatorname{mod}2$)
\end{center}

\noindent and%

\[
\lim_{n\rightarrow\infty}\dfrac{_{L}\mu_{n}}{n}=\,_{L}G_{1/2}%
(1,1)=0.75782301126849283774...,
\]%
\[
\lim_{n\rightarrow\infty}\dfrac{_{L}\sigma_{n}^{2}}{n^{2}}=\,_{L}%
G_{1/2}(1,2)-\,_{L}G_{1/2}(1,1)^{2}=0.03700721658229030320...,
\]%
\[
\lim_{n\rightarrow\infty}\dfrac{_{L}\nu_{n}}{n}=\frac{4e}{(1+e)^{2}%
}=0.78644773296592741014...,
\]
\pagebreak%
\[%
\begin{array}
[c]{ccc}%
\lim\limits_{n\rightarrow\infty}\dfrac{_{S}\mu_{n}}{n^{1/2}}=\sqrt{2}%
\,_{S}G_{1/2}(1,1)=2.06089224152016653900... &  & \text{for EV,}%
\end{array}
\]%
\[%
\begin{array}
[c]{ccc}%
\lim\limits_{n\rightarrow\infty}\dfrac{_{S}\sigma_{n}^{2}}{n^{3/2}}=\sqrt
{2}\,_{S}G_{1/2}(1,2)=1.40007638550124502818... &  & \text{for EV,}%
\end{array}
\]%
\[
\lim\limits_{n\rightarrow\infty}\dfrac{_{S}\mu_{n}}{n^{1/2}}=\left\{
\begin{array}
[c]{lll}%
0.55... &  & \text{for OD and }n\equiv0\operatorname{mod}2,\\
1.50... &  & \text{for OD and }n\equiv1\operatorname{mod}2;
\end{array}
\right.
\]%
\[
\lim\limits_{n\rightarrow\infty}\dfrac{_{S}\sigma_{n}^{2}}{n^{3/2}}=\left\{
\begin{array}
[c]{lll}%
0.12... &  & \text{for OD and }n\equiv0\operatorname{mod}2,\\
1.27... &  & \text{for OD and }n\equiv1\operatorname{mod}2.
\end{array}
\right.
\]

\noindent As before, no explicit integrals are known for the latter four results.

\section{Square Permutations}

A permutation is a square ($p=q^{2}$) if and only if, for any integer $\ell$,
the number of cycles of length $2\ell$ (in its disjoint cycle decomposition)
must be even \cite{Bn-tcs7, Wf-tcs7, O4-tcs7}. \ There is no restriction on
the number of cycles of length $2\ell+1$. \ Hence both%
\[%
\begin{array}
[c]{ccccc}%
(1\;2\;3\;4)(5\;6\;7\;8)(9) &  & \text{and} &  & (1\;2\;3)(4\;5\;6)(7\;8\;9)
\end{array}
\]
are squares, but%
\[%
\begin{array}
[c]{ccccc}%
(1\;2)(3\;4)(5\;6)(7\;8\;9) &  & \text{and} &  & (1)(2\;3\;4\;5)(6\;7\;8)(9)
\end{array}
\]
are not squares. \ While it is not possible to enforce restrictions on
$(2\ell)$-cycle counts using the values of $c_{n}$ alone, we can still employ
brute force methods to calculate%
\[%
\begin{array}
[c]{ccc}%
\{L_{k,3}\}_{k=1}^{3}=\{1,0,2\}, &  & \{S_{k,3}\}_{k=1}^{3}=\{1,0,2\}
\end{array}
\]%
\[%
\begin{array}
[c]{ccc}%
\{L_{k,4}\}_{k=1}^{3}=\{1,3,8\}, &  & \{S_{k,4}\}_{k=1}^{2}=\{9,3\};
\end{array}
\]%
\[%
\begin{array}
[c]{ccc}%
\{L_{k,5}\}_{k=1}^{5}=\{1,15,20,0,24\}, &  & \{S_{k,5}\}_{k=1}^{5}%
=\{36,0,0,0,24\};
\end{array}
\]%
\[%
\begin{array}
[c]{ccc}%
\{L_{k,6}\}_{k=1}^{5}=\{1,45,80,0,144\}, &  & \{S_{k,6}\}_{k=1}^{3}%
=\{230,0,40\};
\end{array}
\]%
\[%
\begin{array}
[c]{ccc}%
\{L_{k,7}\}_{k=1}^{7}=\{1,105,560,0,504,0,720\}, &  & \{S_{k,7}\}_{k=1}%
^{7}=\{960,210,0,0,0,0,720\}.
\end{array}
\]
A significantly faster algorithm might provide insight on the exp-log $a$
parameter (numerical bounds, if relevant) and corresponding statistics.

\section{Addendum}

We justify two of the limiting median $_{L}\widetilde{\nu}_{n}$ formulas.
\ For $a=2$,%
\[%
{\displaystyle\int\limits_{x}^{1}}
\frac{1-y}{y}dy=\left.  \ln(y)-y\right\vert _{x}^{1}=\frac{1}{4}%
\]
when $-\ln(x)+x=5/4$, i.e., $x^{-1}e^{x}=e^{5/4}$, i.e., $(-x)e^{-x}%
=-e^{-5/4}$. \ For $a=3/2,$%
\[%
{\displaystyle\int\limits_{x}^{1}}
\frac{\sqrt{1-y}}{y}dy=%
{\displaystyle\int\limits_{0}^{1-x}}
\frac{\sqrt{z}}{1-z}dz=\left.  -2\sqrt{z}+2\operatorname{arctanh}\left(
\sqrt{z}\right)  \right\vert _{0}^{1-x}=\frac{1}{3}%
\]
when $x=\operatorname{sech}(\xi)^{2}$; this is true because%
\[
-\sqrt{1-x}+\operatorname{arctanh}\left(  \sqrt{1-x}\right)  =-\tanh
(\xi)+\operatorname{arctanh}\left(  \tanh(\xi)\right)  =-\tanh(\xi)+\xi
=\frac{1}{6}.
\]

\section{Acknowledgements}

I am grateful to Mike Spivey \cite{Sp-tcs7}, Alois Heinz and\ Michael Somos
for helpful discussions. \ The creators of Mathematica, as well as
administrators of the MIT Engaging Cluster, earn my gratitude every day.
\ A\ sequel to this paper will be released soon \cite{Fn-tcs7}.

\end{document}